\documentstyle[leqno,amsfonts,12pt]{article}

\newcommand{\CC}{\Bbb C}
\newcommand{\rbox}{$\:\:$ \raisebox{-1ex}{$\:\Box\:$}}

\newcommand{\mb}{\mbox}

\newcommand{\beq}{\begin{equation}}
\newcommand{\eeq}{\end{equation}}

\newcommand{\ve}{\varepsilon}

\newcommand{\ov}{\overline}

\newcommand{\ra}{\rightarrow}

\newcommand{\z}{\zeta}

\newcommand{\p}{\preceq}

\newtheorem{th}{Theorem}

\newtheorem{lem}{Lemma}

\def\ve{\varepsilon }

\newcommand{\ueberschrift}{\bigskip\goodbreak\noindent\bigskip}
\newcounter{theabsatz}
\newcommand{\absatz}[1]{\stepcounter{theabsatz} \ueberschrift
               {\large \bf \arabic{theabsatz}. {#1}} \setcounter{equation}{0}}

\title{On zeros of polynomials orthogonal over a convex domain
\thanks{
This research was supported in part by the National Science
Foundation grant DMS-9707359.} }

\author{V. V. Andrievskii, I. E. Pritsker and R. S. Varga}

\date{}

\begin{document}
\mathsurround=2pt

\maketitle

\begin{abstract}
We establish a discrepancy theorem for signed measures, with a
given positive part, which are supported on an arbitrary convex
curve. As a main application, we obtain a result concerning the
distribution of zeros of polynomials orthogonal on a convex
domain.
\end{abstract}

\absatz{Introduction and main results}

Let $G\subset \Bbb C$ be a bounded Jordan domain, and let $h(z)$
be a weight function on $G $, i.e., a function, which is positive
and measurable on $G$. Next, let $Q_n(z) = Q_n(h,z) = \lambda_n
z^n + \ldots,\, \lambda_n > 0$, $n = 0,1,\ldots$, be the sequence
of polynomials orthogonal in $G$ with respect to the weight
function $h(z)$, that is,
   $$
   \int\limits_G Q_k(z) \overline{Q_l(z)}\, h(z)\, dm(z) =
   \left\{
   \begin{array}{lcl} 1, & \mb{if }\; k=l, \\
                      0, & \mb{if }\; k\neq l,
    \end{array} \right.
   $$
where $dm(z)$ denotes 2-dimensional Lebesgue measure (area).

With $L$ denoting the boundary of $G$, we assume  that
   \beq
   \label{GL_3_1}
   h(z) \ge c\, (\mb{dist}(z,L))^m,  \quad z \in G,
   \eeq
for some constants $m>0,\, c > 0$.

Recently, Eiermann and Stahl \cite{eiesta} made computations and
raised some conjectures about the distribution of the zeros of the
orthogonal polynomials $\tilde{Q}_n(z):=Q_n(h,z)$, in the special
case where $h(z)\equiv 1$, on convex domains $G$ having polygonal
boundaries. In particular, $N$-gons $G_N,\, N=3,4,\ldots ,$ which
have their vertices at the $N$-th roots of unity, were also
considered in \cite{eiesta}. It was previously shown in
\cite{andbla1} that for {\it some} $G$ and {\it some} $n$, the
distribution, of zeros of the associated orthogonal polynomials
$\tilde Q_n$, is governed by the equilibrium measure $\mu_{\ov G}$
of $\ov G$. The main purpose of this paper is to prove a
discrepancy theorem for a special measure $\tau_n$, which is
closely connected with zeros of $Q_n$ and $\mu_{\ov G}$, for {\it
all} convex domains $G$ and $n \in {\Bbb N}$.

In what follows, we assume that $G\subset \Bbb C$ is always
convex. It is known (cf. Stahl and Totik \cite[p. 31]{statot})
that the zeros $ z_{n,1},\ldots ,z_{n,n}$ of $Q_n$ belong to $\ov
G$, for any $n \in {\Bbb N}$.

Let $\omega (z,J,G)$, $z\in G$ and $J\subset L:=\partial G$ be the
harmonic measure of $J$ at $z$ with respect to $G$. We extend this
notion to the boundary points $z\in L$, by setting $$ \omega
(z,J,G):=\left\{\begin{array}{ll} 1,&\quad z\in J,\\[2ex] 0,&\quad
z\not\in J.
\end{array}\right.
$$ Next, we associate with $Q_n$ the measure $$
\tau_n(J):=\frac{1}{n}\sum\limits_{j=1}^n\omega (z_{n,j},J,G),
\quad n \in {\Bbb N}. $$

We will compare $\tau_n$ with the  equilibrium measure $\mu =
\mu_{\ov G}$ of $\ov G$ (see \cite{saftot}), which has a simple
interpretation using the conformal mapping $\Phi $ of $\Omega
:=\overline{\Bbb C}\setminus \ov G$ onto $\Delta :=\{w:|w|>1\}$,
normalized by the conditions $$
   \Phi(\infty) = \infty \quad \mbox{and} \quad \Phi'(\infty)
   := \lim_{z\ra \infty} \;   \frac{\Phi(z)}{z} > 0,
$$ where we define $\Psi \glossary{$\Psi$} :=\Phi ^{-1}. $ Namely,
$\Phi $ can be extended to a homeomorphism $\Phi
:\overline{\Omega} \rightarrow \overline{\Delta} $ and, for any
subarc $J\subset L$, $$ \mu (J)=\frac {1} {2\pi }| \Phi (J)|, $$
where $ |\gamma |\glossary{$|\gamma |$}$ denotes the length of
$\gamma \subset \Bbb C$.

\noindent{\bf Remark.} It is known that the measures $\tau_n$
converge to $\mu_{\ov G}$ in the weak* topology, as $n \to
\infty,$ for any Jordan domain $G$ (cf. Theorem 2.2.1 of \cite[p.
42]{statot} and its proof).

We define the  discrepancy of a signed (Borel) measure $\sigma$,
supported on $L$, by $$ D[\sigma] \glossary{$D[\sigma]$}
:=\sup |\sigma (J)|,
$$ where the supremum is taken over all subarcs $J\subset L$. With
this definition, our new result, for the asymptotic zero
distribution of polynomials orthogonal over a general convex
domain, is stated as

\begin{th}
\label{th1} Let $G$ be a bounded convex domain, and let $h(z)$
satisfy (\ref{GL_3_1}). Then for each $n=2,3,\ldots$, $$ D[\mu
_{\ov G}-\tau_n]\le c\,\sqrt{\frac{\log n}{n}} $$ for some
constant $c>0$, which is independent of $n$.
\end{th}

The main idea of the proof of Theorem \ref{th1} is in its
potential theoretical interpretation. Namely, let cap$\ov G$ be
the (logarithmic) capacity of $\ov G$. We consider the logarithmic
potentials of  $\mu$ and $\tau_n$ in $\Omega$:
    \begin{eqnarray*}
    U(\mu , z) &:=& - \int \log |z - \zeta|\,
    d\mu(\zeta) \\
    &=&
    - \log |\Phi(z)| - \log (\mb{cap }\ov G),
    \end{eqnarray*}
    \begin{eqnarray*}
    U(\tau_n,z) &:=& - \int \log |z - \zeta|\, d\tau_n(\zeta) \\
    &=&
    - \int \log |z - \zeta|\, d\nu_{Q_n}(\zeta) =
    - \frac{1}{n} \log
    \frac{|Q_n(z)|}{\lambda_n}
    \end{eqnarray*}
(where we have used the fact that $\tau_n$ is the balayage of the
zero-counting measure $\nu_{Q_n}$ which associates the mass $1/n$
with each zero of $Q_n$ according to its multiplicity),
 and their difference
    \begin{eqnarray*}
    U(\mu-\tau_n,z) &:=& U(\mu ,z) - U(\tau_n,z)
    \\
&=&
\frac{1}{n}\log \frac{|Q_n(z)|}{\lambda_n(\mb{cap }\ov G)^n|\Phi (z)|^n}.
\end{eqnarray*}

It is  proved in \cite{andbla1} that the inequalities \beq
\label{n.1} ||Q_n||_{\ov G}:=\sup\limits_{z\in \ov G}|Q_n(z)|\le
c_1\, n^{c_2}, \eeq \beq \label{n.2} \lambda_n\, (\mb{cap }\ov
G)^n\ge c_3\, n^{-2} \eeq hold for some constants $c_j>0,j=1,2,3$,
which are independent of $n$. This implies that, for any $n\ge 2$,
$$ U(\mu -\tau_n,z)\le c_4\,\frac{\log n}{n}\, ,\quad z\in\Omega,
\quad c_4>0, $$ where $c_4$ is also independent of $n$.

Theorem \ref{th1} is actually a consequence of our result given
below, which is a new Erd\H{o}s-Tur\'an-type theorem (its proof
will be given in subsequent sections).

\begin{th}
\label{th2} Let $G\subset\Bbb C$ be a bounded convex domain, and
let $\tau$ be a unit Borel measure supported on $L:=\partial G$.
If $$ \varepsilon =\varepsilon (\tau ):=\sup\limits_{z\in\Omega}
U(\mu_{\ov G}-\tau ,z)\,\, (\ge 0), $$ then \beq \label{n.5}
D[\mu_{\ov G} -\tau ]\le c\,\sqrt{\varepsilon} , \eeq for some
constant $c>0$, independent of $\tau$.
\end{th}

For $G={\Bbb D}:=\{ z:|z| < 1\}$, the result of Theorem \ref{th2}
is due to Ganelius \cite{gan}, which in turn generalized results
of Erd\H{o}s and Tur\'an \cite{erdtur}, concerning distribution of
zeros of polynomials with given uniform norms on the unit disk.
Further results and bibliography of papers devoted to this subject
can be found in \cite{blagro,blamha,tot,andbla,saftot}.

The following example shows the sharpness of Theorem \ref{th2}.

{\bf Example 1}. Let $G=\Bbb D$ and let $\mu_\delta ,0<\delta \le
1$, be the equilibrium measure of $V_\delta:=\ov{\Bbb D}\cup
[1,1+\delta]$. Consider the measure $\tau_\delta$, supported on
the unit circle $\Bbb T:=\partial \Bbb D$, which is defined for
any Borel set $B\subset\Bbb T$ by the formula $$ \tau_\delta
(B):=\mu_\delta(\{ z\in {\Bbb C} \setminus\{ 0\} : z/|z|\in B\} ).
$$ It is easy to see that $$ \mb{cap }V_\delta =\frac 14\left(
3+\delta+\frac{1}{1+\delta}\right) =1+\frac{\delta^2}{4(1+\delta
)}\, . $$ Therefore for $z\in\Bbb T$ we have $$ U(\mu-\tau_\delta
,z)\le U(\mu -\mu_\delta ,z)=\log \mb{cap }V_\delta
\le\frac{\delta^2}{4}\, . $$ At the same time an elementary
computation, involving the transformation $z\to (z+1/z)/2$, shows
that $$ D[\mu -\tau_\delta]\ge |(\mu -\tau_\delta )(1)|=
\mu_\delta ([1,1+\delta ])\ge \frac{\delta}{3\pi}\, . $$ This
implies that $$ D[\mu -\tau_\delta]\ge \frac{2}{3\pi}\sqrt{\ve
(\tau_\delta )}, $$ which shows the sharpness of Theorem
\ref{th2}. \hfill\rbox

Note that statements similar to Theorem \ref{th1} can also be
proved (by making of use of Theorem \ref{th2}) for other systems
of polynomials. All that is needed for this purpose is to
establish the analogues of
 (\ref{n.1}), (\ref{n.2}) and the assumption that
\beq
\label{n.3}
all\,
zeros\,\, of\,\, the\,\, corresponding\,\,
 polynomials\,\, belong\,\, to\,\, \ov G.
\eeq We cite three examples of well-known polynomials suited for
such applications of Theorem \ref{th2}. In all of them, $G$ is a
convex domain and $n \in \Bbb N$.

{\bf Example 2}. Let $F_n(z):=($cap$ \ov G)^{-n}z^n+\ldots $ be
the $n$-th Faber polynomial for $\ov G$ (cf. \cite{smileb}). Then,
(\ref{n.3}) is valid by \cite[Theorem 2]{koepom}. In addition, we
have by the same Theorem 2 of \cite{koepom} that $$ ||F_n||_{\ov
G}\le 2, \quad n \in \Bbb N.$$

{\bf Example 3}. Consider the derivatives $F'_{n+1}(z)$ of the
above Faber polynomials. For these polynomials, condition
(\ref{n.3}) is then proved in \cite{ull}. At the same time, by the
Markov-type inequality for complex polynomials, which is a simple
consequence of L\"{o}wner's distortion theorem (see, for example,
\cite[p. 58]{andbeldzj}), there holds $$ ||F'_{n+1}||_{\ov G}\le
c\,(n+1)^2,\quad c=c(G)>0. $$

{\bf Example 4}. Let $T_n(z)=z^n+\ldots$, be the $n$-th normalized
Chebyshev polynomial for $\ov G$. Condition (\ref{n.3}) is then
well known (cf. \cite{smileb}). The corresponding estimate for the
uniform norm on $\ov G$ follows from the extremal property of
Chebyshev polynomial: $$ ||T_n||_{\ov G}\le (\mb{cap }\ov G)^n
||F_n||_{\ov G} \le 2\ (\mb{cap }\ov G)^n. $$

In what follows, we denote by $c,c_1,\ldots$ positive constants
and by $\ve_0,\varepsilon_1,\ldots$ sufficiently small positive
constants (different each time, in general) that either are
absolute or depend on parameters not essential for the arguments;
sometimes such a dependence will be indicated. For $a > 0$ and $b
> 0$ we use the expression $a \preceq b$ (order inequality) if $a
\le c\, b$ for some $c > 0$. The expression $a \asymp b$ means
that $a \preceq b$ and $b \preceq a$ hold simultaneously.

\absatz{Some facts from geometric function theory}

Each convex curve is known to be quasiconformal (see \cite[pp. 63,
87]{leh}). It is further known (see \cite[Chapter IV]{ahl}) that
the conformal mapping $\Phi$ can be extended in this case to a
quasiconformal mapping of the whole plane onto itself. We keep the
same notation for this extension. Note that the inverse function
$\Psi
:=\Phi^{-1}$
 will be quasiconformal too.

The following result is useful in the study of metric properties
of the mappings $\Phi$ and $ \Psi$.
%*****************
\begin{lem}
\label{q} (\cite[p. 97]{andbeldzj}) Let $w=F(\zeta )$ be a
$K$-quasiconformal mapping of $\ov{\Bbb C}$ onto itself with
$F(\infty )=\infty ,\, \zeta_j\in {\CC},\, w_j:=F(\zeta_j),\,
j=1,2,3$, and $|w_1-w_2|\le c_1|w_1-w_3|$. Then $|\zeta _1-\zeta
_2|\le c_2\, |\zeta_1-\zeta _3|$ and, in addition, $$\left| \frac
{\zeta _1-\zeta _3}{\zeta _1-\zeta _2}\right|\le c_3\left|
\frac{w_1-w_3}{w_1-w_2}\right|^K,$$ where $c_j=c_j(c_1,K),\,
j=2,3$.
\end{lem}

The convexity of $G$ implies some special distortion properties of
the function $\Phi$.
%*********************
\begin{lem}
\label{a} Let $z_1 \in L$, $z_2,z_3\in\ov{\Omega}$ and
$\tau_j:=\Phi (z_j),j=1,2,3$. If $|\tau_1 - \tau_2| \le |\tau_1 -
\tau_3| \le 1$, then the inequality
   \beq
   \label{3.1}
\left|\frac{z_2-z_1}
   {z_3-z_1}\right|
   \le c_1
   \left|\frac{\tau_2-\tau_1}{\tau_3-\tau_1}\right|
   \eeq
holds with $c_1 = c_1(G) > 0$.
\end{lem}
%****************

{\bf Proof.} Without loss of generality, we assume that
    $$
     |z_2-z_1| < |z_3-z_1| < \frac{1}{2}\mb{ diam }L
    $$
(otherwise (\ref{3.1}) follows easily from Lemma \ref{q}). Next we
introduce the following notations. Denote by $\gamma(x)=\gamma
(z_1,x) \subset \Omega$, for $0 < x < \frac{1}{2}\mb{ diam }L$,
the subarc of the circle $\{\xi:|\xi-z_1| = x\}$ that separates
the point $z_1$ from $\infty$ in $\Omega$. Let
$Q(\delta,t)=Q(z_1,\delta,t)$, for $0
< \delta < t < \frac{1}{2}\mb{ diam }L$, be the quadrilateral
bounded by the arcs $\gamma (\delta)$, $\gamma (t)$ and the two
subarcs of $L$ joining their endpoints. Denote the family of all
locally rectifiable arcs in $Q(\delta,t)$, which separate the
sides $\gamma (\delta)$ and $\gamma (t)$, by $\Gamma (\delta,t)$,
and the module of $\Gamma (\delta,z)$ by $m(\delta,t)$ (see
\cite{ahl,lehvir}). By the comparison principle $$
   m(\delta,t) \le \frac{1}{\pi} \log \frac{t}{\delta} , \qquad
   0 < \delta < t < \frac{1}{2}\mb{ diam }L.
$$ For any triplet of points $\xi_1, \xi_2, \xi_3 \in
\overline{\Omega}$ with $|\xi_1 - \xi_2| = |\xi_1 - \xi_3|$, we
have by Lemma \ref{q} that $$ |\Phi(\xi_1) - \Phi(\xi_2)|\asymp
|\Phi (\xi_1) - \Phi(\xi_3)|. $$ Hence, according to \cite{belmik}
(see also \cite[p. 36]{andbeldzj})
   $$
   \left|\frac{\tau_3- \tau_1}{\tau_2 - \tau_1}\right| \asymp \exp
   (\pi m (|z_2-z_1|, |z_3-z_1|))
   \le
   \left|\frac{z_3 - z_1}{z_2 - z_1}\right|.
   $$

\hfill\rbox

\begin{lem}
\label{lem1.12} The inequality \beq \label{GLA_1_51} \omega
(z,l,{\Bbb D})\le 8\,\frac {1-|z|}{\mb{\em dist }(z,l)} \eeq holds
true for any $z\in \Bbb D$ and any arc $l\subset \Bbb T$.
\end{lem}

{\bf Proof}. Using a rotation with respect to the origin, we can
reduce the situation to the case when $0<z<1$ and $l=\{
e^{i\theta} : \theta_1\le\theta\le\theta_2\}$,
$0<\theta_1<\theta_2<2\pi +\theta_1$. Moreover, we can assume that
$\theta_2<2\pi$ (since in the other case (\ref{GLA_1_51}) is
trivially valid). Set $$ l_1:=\{ \z\in l:\mb{ Im }\z \ge 0\} ,
\quad l_2:=l\setminus l_1. $$ We assume that $l_1\neq\emptyset$. A
simple geometric reasoning shows that, for $\z=e^{i\theta} \in
l_1$, $$ |\z -z|\ge\frac 1 \pi (\theta -\theta_1),\quad |\z -z|\ge
|z-z_1|,\quad z_1:=e^{i\theta_1}. $$ Therefore, by the Poisson
formula
\begin{eqnarray*}
\omega (z,l_1,\Bbb D)&=&
\frac{1-|z|^2}{2\pi}\int\limits_{l_1}\frac{|d\z |}{|\z -z|^2}
\le\frac{1}{\pi}(1-|z|)\int\limits_{l_1}\frac{|d\z |}{|\z -z|^2}\\
&\le&
4\pi (1-|z|)\int\limits_{\theta_1}^{\pi}
\frac{d\theta}{(\pi |z-z_1|+\theta -\theta_1)^2}\\
&\le&
\frac{4(1-|z|)}{|z-z_1|}\le
\frac{4(1-|z|)}{\mb{dist }(z,l)}.
\end{eqnarray*}
Writing the same estimate for $\omega (z,l_2,\Bbb D)$ and taking
their sum, we get (\ref{GLA_1_51}). \hfill\rbox

\absatz{Auxiliary results}

In this section, we discuss the results needed in the proof of
Theorem 2.

The concept of a regularized distance to an arbitrary compact set
$E\subset {\Bbb R}^n$ is described in \cite[pp. 170-171]{ste}. It
is based on the decomposition of open sets into cubes and the
partition of  unity, which is due to Whitney. It is enough for our
purposes to assume that $E$ is a continuum in the complex plane,
with the simply connected complement $U$. In this case, the notion
of a regularized distance can be explained by making use of the
properties  of a conformal  mapping of $U$ onto the  unit disc.

Namely, let $U\subset {\Bbb C}$ be a simply connected domain,
$E:=\ov{\Bbb C}\setminus U\neq\emptyset,$ with $\infty \in E$.
Denote the distance from $z$ to $E$ by $d(z):=d(z,E)$. This
function is in general not smoother on $U$ than what the obvious
Lipschitz-condition-inequality $$ |d(z)-d(\z )|\le |z-\z |,\quad
z,\z\in\Bbb C, $$ indicates.

It is desirable for several applications to replace $d(z)$ by a
regularized distance $\rho (z)$, which is infinitely
differentiable for $z\in U$. In addition, this regularized
distance should have essentially the same behavior as $d(z)$.

Let $g:U\to H_+:=\{w: \mb{Im\,}w>0\}$ be a conformal mapping. Set
$u(z):=\mb{Im\,} g(z)$. The function \beq \label{a4.0} \rho
(z):=\frac{u(z)}{|g'(z)|},\quad z\in U, \eeq is called a {\it
regularized distance} from $z$ to $E$.

\begin{lem}
\label{lem1.3} (\cite[Lemma 1]{ad}) For each $z\in U$, we have
\beq \label{GLA_1_5} \frac {1}{4}\, \frac {u(z)}{d(z)}\le |g'
(z)|\le 4\, \frac {u(z)}{d(z)}. \eeq Moreover, if $|\xi -z|\le
d(z)/2$ then \beq \label{GLA_1_6} \frac {1}{16}\, \frac
{u(z)}{d(z)} |\xi -z|\le |g(\xi ) -g(z)|\le 16 \, \frac
{u(z)}{d(z)} |\xi -z|. \eeq
\end{lem}

Applying (\ref{GLA_1_5}) we have
$$
\frac{1}{4}\, d(z)\le\rho (z)\le 4d(z),\quad
z\in U.
$$

We note the following fact about the smoothness properties of
$\rho (z)$. Let $f(z),\, z=x+iy$, be a non-vanishing analytic
function in $U$. A simple calculation yields that, for any $z\in
U$,
\beq \label{a4.2} |f|'_x=|f|\, (\log |f|)'_x=|f| \mb{ Re
}(\log f)_z'=|f|\mb{ Re }\frac{f'_z}{f}\, , \eeq \beq \label{a4.3}
|f|'_y=|f|\, (\log |f|)'_y=|f| \mb{ Re }(i\log f)_z'=-|f|\mb{ Im
}\frac{f'_z}{f}\, ; \eeq whence, we conclude that \beq \label{a4.4}
\left| |f|_\xi'\right|\le |f_z'|\quad (\mb{with } \xi =x\,\mb{ or
}\xi=y). \eeq Formulae (\ref{a4.2}) and (\ref{a4.3}) imply that
$\rho (z)\in C^\infty (U)$. Differentiating them once more, we
obtain for $j+k=2,\, j,k\ge 0$, that \beq \label{a4.5}
\left|\frac{\partial^2|f|}{\partial x^j\partial y^k}\right| \le
|f''_{zz}|+2\frac{|f'_z|^2}{|f|}\, . \eeq

Next, we claim that for $z=x+iy\in U;\,
j,k=0,1,2,\, 1\le j+k\le 2$,
\beq
\label{a4.6}
\left|\frac{\partial^{j+k}}{\partial x^j\partial y^k}\,
\rho (z)\right|\le c_1\,\rho (z)^{1-j-k}
\eeq
for some absolute constant $c_1>0$.

Indeed, inequality (\ref{a4.6}) follows immediately from
(\ref{a4.4}), (\ref{a4.5}) and (\ref{GLA_1_5}) after a twice
repeated  differentiation of the formula (\ref{a4.0})
 with respect to $\xi_j=x$ or $\xi_j=y$, $j=1,2$: $$
\frac{\partial\rho}{\partial\xi_1}=\frac{1}{|g'_z|^2}
\left(u'_{\xi_1}|g'_z|-u|g'_z|'_{\xi_1}\right) , $$
\begin{eqnarray*}
\frac{\partial^2\rho}{\partial\xi_1\partial\xi_2}
&=&\frac{1}{|g'_z|^4}\left\{\left( u''_{\xi_1\xi_2}|g
'_z|+u'_{\xi_1}|g'_z|'_{\xi_2} \right.\right.\\ &-&\left.
u'_{\xi_2}\, |g'_z|'_{\xi_1}-u\, |g'_z|''_{\xi_1\xi_2}\right)
|g'_z|^2 \\ &-& \left. 2\left( u'_{\xi_1}\, |g'_z|-u\,
|g'_z|'_{\xi_1}\right) |g'_z|\, |g'_z|'_{\xi_2}\right\}\, ,
\end{eqnarray*}
if we know that for $k=2,3$, \beq \label{a4.7} |g^{(k)}(z)|\le
c_2\,u(z)\,\rho (z)^{-k},\quad z\in U, \eeq with an absolute
constant $c_2>0$.

In order to prove (\ref{a4.7}), we put $d:=d(z)/32$ and note that
by (\ref{GLA_1_6}),  $$ |g (\z )-g (z)|\le\frac{1}{2}\, u(z), $$
for any $\z$ with $|\z -z|=d$. Therefore, we have, according to
(\ref{GLA_1_5}), that $$ |g'(\z )|\le 4\,\frac{u(\z )}{d(\z )} \le
10\, \frac{u(z )}{d(z )}\, , $$ for such $\z$. Next, we apply
Cauchy's formula and (\ref{GLA_1_5}) to obtain that for $k=2,3$,
\begin{eqnarray*}
|g^{(k)}(z)|&=&\frac{(k-1)!}{2\pi}\left|\, \int\limits_{|\z
-z|=d}\frac{g'(\z )}{(\z -z)^k}\: d\z\,\right|\\ &\le& 10\,
(k-1)!\,\, 32^{k-1}\, \frac{u(z)}{d^k(z)}\, .
\end{eqnarray*}
This completes the proof of (\ref{a4.7}) and, consequently, of (\ref{a4.6}).

The second topic concerns the ``body-contour" properties of
harmonic functions. Let $G\subset\Bbb C$ be a bounded convex
domain, and let $f(z)$ be a real valued function, which is
continuous on $\ov{G}$ and harmonic in $G$. Let $z\in L:=\partial
G, \z\in G,\delta :=|z-\z|$. We next estimate the quantity
$|f(\z )-f(z)|$ in terms of the local modulus of continuity of $f$
on $L$, that is, $$ \omega_{z,f,L}(t):=\sup\limits_{\xi\in L\atop
|\xi -z|\le t} |f(\xi )-f(z)|,\quad t>0. $$ Let $z_0\in G$ be a
fixed point. We assume that $2\delta<$ dist $(z_0,L) =:d_0$. For
$0 <t<d_0$, denote by $\gamma (t)=\gamma (z,t)$ a crosscut of $G$,
i.e. an open Jordan arc in $G$ with endpoints on $L$, which is a
subarc of the circle $\{\xi :|\xi -z|=t\}$ and has nonempty
intersection with the interval $[z,z_0]$. The endpoints of $\gamma
(t)$ divide $L$ into two subarcs. Denote the subarc containing $z$
by $l(t)$.

Since $L$ is quasiconformal, Ahlfors' geometric criterion (see
\cite{ahl}) gives the inequality
   \beq
   \label{n.10}
   \min\{\mb{diam } L', \mb{diam } L''\} \le c\, |z_1 - z_2|,
   \quad \mb{ for any } z_1,z_2 \in L,
   \eeq
with $c = c\, (L) \ge 1$, where $L'$ and $L''$ are the associated
two arcs $L \backslash \{z_1,z_2\}$ consists of. Therefore, the quantity $$
M=M(z_0,L):=\sup\limits_{z\in L}\sup\limits_{0<t<d_0}
\frac{\mb{diam }l(t)}{t} $$ is finite. Moreover, it is easy to
prove that $M\le M_0$, where $M_0$ depends only on the constant
$c$ from (\ref{n.10}), and consequently only on the constant of
quasiconformality of $L$.

Let $$ \nu (t):=\omega (\z ,L\setminus l(t),G),\quad 0 <t
<d_0, $$ be the corresponding harmonic measure. Next, we fix a
number $s$, satisfying $2\delta <s\le d_0$, and define a natural
number $k$ such that $$ \frac{s}{2}\le 2^k\,\delta <s. $$ By the
maximum principle for harmonic functions, we have
\begin{eqnarray*}
|f(\z )-f(z)| &\le & \omega_{z,f,L}(M\delta ) +\sum_{j=0}^{k-1}
\omega_{z,f,L}(M2^{j+1}\,\delta )\,\nu (2^j\,\delta ) +2\,
||f||_L\,\nu \left( \frac s2\right)\\ &\le& \omega_{z,f,L}(M\delta
)+2\int\limits_\delta^s \frac{\omega_{z,f,L}(2Mt)}{t}\,\nu
\left(\frac t2\right) dt +2||f||_L\,\nu \left( \frac s2\right) .
\end{eqnarray*}

Our next goal is to obtain effective estimates of the harmonic
measure $\nu (t)$. Let $\Gamma =\Gamma (\z ,l(t),G), \delta
<t<d_0$, be a family of all crosscuts of $G$ that separate point
$\z$ from $L\setminus l(t)$. We note that \beq \label{eq3.11}
m(\Gamma )\le\frac{1}{\pi}\log\frac{4}{\nu (t)}. \eeq

Indeed, taking into account that both module and harmonic measure
are conformal invariants, we introduce the conformal mapping
$g:G\to \Bbb D$ such that $$ g(\z )=0,\quad g(L\setminus l(t))=\{
e^{i\theta\pi}:-a\le \theta\le a\},\quad a:=\nu (t). $$

According to \cite[pp. 319--320]{her} (see also \cite[p. 6]{hal}),
we have $$ m(\Gamma)^{-1}=m(g(\Gamma))^{-1}= 2T\left(\sin
\frac{\pi}{2}(1-a)\right) =2T\left(\cos \frac{\pi a}{2}\right) ,
$$ where we set $$ T(k):=\frac{K((1-k^2)^{1/2})}{K(k)} $$ and $$
K(k):=\int\limits_0^1(1-x^2)^{-1/2}(1-k^2x^2)^{-1/2}dx,$$ for
$0<k<1$. Hence $$ 2m(\Gamma )=T\left(\sin \frac{\pi a}{2}\right) .
$$ By \cite[p. 61]{lehvir}, $$ T\left(\sin \frac{\pi a}{2}\right)
\le \frac{2}{\pi}\log \frac{4}{\sin \frac{\pi a}{2}} \le
\frac{2}{\pi}\log \frac{4}{a}. $$ Thus, we obtain (\ref{eq3.11})
by comparing the last two equations.

On the other hand, comparing the families $\Gamma$ and
$\Gamma_1:=\{\gamma (u)\}_{\delta <u<t}$, we have $$ m(\Gamma )\ge
m(\Gamma_1)\ge \frac{1}{\pi}\log\frac{t}{\delta}\, . $$ Therefore,
it follows from (\ref{eq3.11}) that  $$ \nu (t)\le 4\,
\frac{\delta}{t}, $$ and that
\begin{eqnarray}
 \nonumber |f(\z )-f(z)| && \\ \label{n3.11} &\le& \omega_{z,f,L}(M\delta
)+16\, \delta\int\limits_\delta^s
\frac{\omega_{z,f,L}(2Mt)}{t^2}dt +16\, ||f||_L\,\frac \delta s\,
.
\end{eqnarray}

\absatz{Proof of Theorem 2}

Let $\sigma :=\mu -\tau$. We can assume that $0 < \ve <\ve_0$, where
$\ve_0=\ve_0(G)$ is small enough for our constructions below. Let
$J \subset L$ be an arbitrary subarc. In order to prove
(\ref{n.5}), it is sufficient to show that
   \beq
   \label{4.2}
   -\sigma (J)
   \le c\, \sqrt{\ve},
   \eeq
 for
$J$ small enough.

We
set
   $$
   \gamma := \Phi(J) = \{e^{i\theta}:\theta_1 \le \theta \le \theta_2\},
   $$
   $$
   \gamma(r) := \{e^{i\theta}: \theta_1 - r \le \theta \le \theta_2
   + r\}, \quad r > 0,
   $$
   $$
   J(r) := \Psi(\gamma(r)), \quad r > 0.
   $$

Next, we introduce a curvilinear sector based on $J$.
 Let $z_0 \in G$ be a fixed point.
 Denote by $w = \varphi(z)$ the conformal mapping of $G$ onto $\Bbb D$ with the
 normalization $\varphi(z_0) = 0$, $\varphi'(z_0) > 0$. Set
  $\psi := \varphi^{-1}$.
Since $L$ is quasiconformal, the  functions $\varphi$ and $\psi$
can be  extended to the quasiconformal mappings of the extended
complex plane $\ov{\Bbb C}$ onto itself with $\infty$ as a fixed
point (see \cite[Chapter IV]{ahl}), where we keep the same
notations for these extensions.

Letting
$$
\varphi(J) = \{e^{i\theta}:\tilde{\theta_1} \le \theta \le
\tilde{\theta_2}\},
$$
we set
\begin{eqnarray*}
B(J)&:=&
\{ \z\in\ov{\Omega}:\theta_1\le \mb{arg }\Phi (\z )\le \theta_2\}
\\
&\cup&
\{ \z\in\ov{G}:\tilde{\theta_1}\le \mb{arg }\varphi (\z )\le
\tilde{\theta_2}\}.
\end{eqnarray*}

Set $t:=\sqrt{\ve}$ and consider the function
   $$
   h(z) := \left\{ \begin{array}{ll}
   1, & \mb{if } z \in B(J(t)),\\[2ex]
   0, & \mb{otherwise in }
   \Bbb C. \end{array} \right.
   $$
Let $\rho (z)=\rho (z,B(J)),\ z \in \Bbb C$, be a  regularized
distance to $B(J)$ (see Section 3), i.e., a function with the
following properties:
   \beq
   \label{4.4}
    \frac{1}{4}\, \mb{dist} (z,B(J)) \le \rho (z) \le 4\,
     \mb{dist} (z,B(J)),
    \quad z \in \Bbb C,
    \eeq
    \beq
    \label{4.5}
    \rho (z) \in C^\infty (\Bbb C),
    \eeq
   \beq
   \label{4.6}
   \left|\frac{\partial^{j+k}}{\partial x^j\partial y^k}\,
    \rho (x+iy)\right| \le
   c\, \rho (x+iy)^{1-j-k}, \quad j+k=1,2.
   \eeq

Next, we average the function $h$ in the following way
   $$
   \displaystyle
   g(z) :=\left\{ \begin{array}{ll} \displaystyle
   \frac{64}{\rho (z)^2} \int_{\Bbb C} h(\z) V\left(\frac{8\, (\z-z)}
   {\rho (z)}\right) dm(\z), & \mb{if }z \in {\Bbb C}
   \setminus B(J),\\[2ex]
   1\, ,& \mb{if }z\in B(J),
   \end{array} \right.
   $$
where $V(\z)$ is an arbitrary symmetric averaging kernel,
 i.e., $V(z)
 \in C^{\infty}(\CC)$,
    $$
    V(z) = V(|z|) \ge 0, \quad z \in \CC ,
    $$
    $$
    V(z) = 0, \quad |z| \ge 1,
    $$
    $$
    \int V(z)\, dm(z) = 1.
    $$

Note that $g \in C^\infty(\Bbb C)$ by virtue of (\ref{4.5}). Set
$$ L_\ve :=\{ z\in \Omega
:|\Phi (z)|=1+\ve\} , $$ $$ z_L:=\Psi (\Phi (z)/|\Phi (z)|),\quad
z\in\Omega\setminus\{\infty\} . $$ By Lemma \ref{q}, there exists
a sufficiently small constant  $\ve _1>0$ such that  $$
\mb{dist}(z,B(J))< \mb{dist}(z,{\Bbb C} \setminus B (J(t)), $$ for
$z\in L_\ve$, with $z_L\in J(2\ve_1t)$. Therefore,  $$
g(z)=1,\qquad  z\in L_\ve,\, z_L\in J(2\ve_1t), $$ according to
(\ref{4.4}). Further, by the same Lemma \ref{q}, there exists a
sufficiently large constant $c_1>0$ such that $$
\mb{dist}(z,B(J))\le 2\, \mb{dist}(z,B(J(t))), $$ for $z\in L_\ve$
with $ z_L\in L\setminus J(c_1t)$. Therefore, we have for such $z$
that $$ \rho (z)\le 4\, \mb{dist}(z,B(J))\le\, 8\,
\mb{dist}(z,B(J(t))), $$ by (\ref{4.4}), and we obtain $$ g(z)
=0,\quad z\in L_\ve ,  \, z_L\in L\setminus J(c_1t). $$

If $z=x+iy$ and $ \xi =\tilde x+i\tilde y\in L_\ve$ with $z_L,\,
\xi_L \in L(\z_3,\z_1)$, where $\z_1:=\Psi(e^{i\theta_1}),\
\z_3:=\Psi (e^{i(\theta_1-3c_1t)})$ and $L(\z_3,\z_1):=\{\z=\Psi
(e^{i\theta}):\ \theta_1-3c_1t \leq \theta \leq \theta_1 \}$, then
we obtain by Taylor's formula that \beq \label{4.10} g(z )=g(\xi
)+A(\xi )(x-\tilde x)+B(\xi )(y-\tilde y)+r(z,\xi ), \eeq where we
have \beq \label{4.11} |A(\xi )|+|B(\xi )|\preceq |\z_1-\z_3|^{-1}
\eeq and \beq \label{4.12} |r(z,\xi )|\preceq \frac{|z-\xi
|^2}{|\z_1-\z_3|^2}\, , \eeq according to (\ref{4.6}).

The same relations are valid for $z,\xi \in L_\ve $ with $z_L,\,
\xi_L \in L(\z_2,\z_4)$, where $\z_2:=\Psi(e^{i\theta_2}),\
\z_4:=\Psi (e^{i(\theta_2+3c_1t)})$ and $L(\z_2,\z_4):=\{\z=\Psi
(e^{i\theta}):\ \theta_2 \leq \theta \leq \theta_2+3c_1t \}$.

We denote  the harmonic extension of $g$ from $L_\ve$ to
$\overline{\Bbb C} \setminus L_\ve$ by $f(z)$. Set $$
\tilde{f}(w):= f(\Psi (w)),\quad w\in \ov{\Delta}. $$ Then the
following estimate holds.

%*****************
\begin{lem}
\label{m}
 Let $1\le |w|\le 1+2\ve$.
Then
    \beq
    \label{4.13}
    |\tilde{f}(w)-\tilde{f}(w_\ve)| \le c_2\, t,\quad
    w_\ve :=\frac{w}{|w|}(1+\ve).
    \eeq
\end{lem}

The proof of Lemma \ref{m} will be given in the next section.

Further, we average the function $\tilde{f}$ in the following
 way.
 Let $V(z)$, $z \in \CC$, be an   averaging kernel as above.
Consider the function
   $$
   \tilde{u}(w) := \left\{ \begin{array}{ll}
   \displaystyle
   \frac{16}{\ve^2} \int
   \tilde{f}(t) V \left( \frac{4(t-w)}{\ve} \right) \, dm(t), &
   \; \mbox{if }
     1 + \frac{3}{4}\ve \le |w| \le 1 + \frac{5}{4} \ve ,\\[2ex]
    \tilde{f}(w), & \; \mbox{elsewhere in } \ov{\Delta}.
    \end{array} \right.
    $$
 Note that $\tilde{u} \in C^\infty (\Delta)$,
     \beq
     \label{4.15}
    0 \le \tilde{u}(w) \le 1, \quad w \in \ov{\Delta},
    \eeq
    and that the Laplacian of $\tilde u$ satisfies
    \beq
    \label{4.16}
    |\Delta \tilde{u}(w)| \preceq \frac{t}{\ve^2} ,\quad
    1 + \frac{3}{4}\ve \le |w| \le 1 + \frac{5}{4} \ve ,
    \eeq
    by (\ref{4.13}).
Let us introduce the function
    $$
    u(z) := \left\{ \begin{array}{cc}
    \tilde{u}(\Phi(z)), & \; \mb{if } z \in \Omega ,\\[2ex]
    f(z), & \; \mb{if } z \in \ov{G},
    \end{array} \right.
    $$
which obviously belongs to the class $ C^{\infty}(\CC)$. It
follows that
   \beq
   \label{4.17}
   \int \Delta u(z)\,  dm(z) = 0,
   \eeq
by Green's formula. Applying the techniques of \cite{blagro}, we
can establish the inequality
   \beq
   \label{4.18}
   \left|\int u\, d\sigma \right| \preceq t .
   \eeq
Indeed, on setting
     $$
    \tilde{U}(\sigma ,w) :=
    U(\sigma , \Psi(w)), \quad
     w\in \Delta ,
     $$
 and, using the representation of the function $u$ by means of Green's formula
    $$
    u(z) = u(\infty) + \frac{1}{2\pi} \int \Delta u(\zeta) \log
    |z-\zeta|\, dm(\zeta ), \quad z \in \CC ,
    $$
 we obtain that
    \begin{eqnarray*}
    \left| \int f\, d\sigma \right| &=& \frac{1}{2\pi}\left| \int
    \left( \ve
    -U (\sigma ,\zeta)\right) \Delta u(\zeta)\, dm(\zeta )\right|
    \\
    &\le& \frac{1}{2\pi} \int \left( \ve -
    \tilde{U}(\sigma ,w)\right) |\Delta \tilde{u}(w)|\, dm(w)
    \p t,
    \end{eqnarray*}
by (\ref{4.16}) and (\ref{4.17}) (see \cite{blagro} for details).

The equations (\ref{4.13}), (\ref{4.15})  and (\ref{4.18}) imply
that
  \begin{eqnarray*}
   -\sigma (J)
   &\le&- \int u  d\sigma + \mu
   (J(c_1t)\setminus J) + \int\limits_{L \backslash J (c_1t )} u d\mu
   \\
   &+ &\int\limits_{J} (1-u) d\tau
   \preceq t\, ,
   \end{eqnarray*}
which is the assertion of (\ref{4.2}).

\hfill\rbox

\absatz{Proof of Lemma \ref{m}}

Let $w=re^{i\theta}$. Applying  Lemma \ref{lem1.12}, we easily
obtain (\ref{4.13}) for the case $$ 1+\ve <r<1+2\ve \, , $$ $$
\theta_1 -\ve_1t\le \theta \le \theta_2+\ve_1 t\quad \mb{or}\quad
\theta_2+2c_1 t\le \theta \le 2\pi + \theta_1-2c_1t\, . $$
If $$ 1+\ve <r< 1+2\, \ve ,\quad \theta_1 -2c_1t\le \theta \le
\theta_1-\ve_1 t\, , $$
 we set
 $ \xi =\z_\ve :=\Psi (w_\ve )$  and write the function $g$
in the form of (\ref{4.10}).

Lemma \ref{q} and Lemma \ref{a} imply that \beq \label{5.2} \left|
\frac{\z-\z_\ve}{\z_1-\z_3}\right| \preceq \left|
\frac{\z_L-\z_\ve}{\z_L-\z_1}\right|\preceq \frac{\ve}{t}=t. \eeq
Define the harmonic extension of the function appearing in
(\ref{4.10}) to ext$L_\ve\setminus\{\infty\}$ by the formula $$
r(z,\xi ):= f(z)-g(\xi )-A(\xi )(x-\tilde x)-B(\xi )(y-\tilde y),
$$ and set $$ \tilde{r}(\tau ):= r(\Psi (\tau ),\xi ), \quad |\tau
|\ge 1+\ve . $$ Note that for $z\in L_\ve$ with $z_L\in
L(\z_3,\z_1)$, we have that \beq \label{n.7} \left| \frac{z-\xi
}{\z_1-\z_3}\right| \preceq \frac{|\Phi (z)-w_\ve |}{t}\, . \eeq
Indeed, without loss of generality we assume that $|z-\z_1|\ge
|z-\z_3|$, and therefore $$ |\z_1-\z_3|\asymp |z-\z_1|, $$ $$ \Phi
(z)=:\tau =(1+\ve )e^{i\eta},\quad |\theta_1-\eta| \asymp t. $$ If
$|\theta -\eta|\ge \ve /32$, then (\ref{n.7}) follows from Lemma
\ref{q} and Lemma \ref{a}, because $$ \left| \frac{z-\xi
}{\z_1-\z_3}\right|\asymp \left| \frac{z_L-\xi
}{z_L-\z_1}\right|\p \left| \frac{\Phi (z_L)-\Phi (\xi )}{\Phi
(z_L)-\Phi (\z_1)}\right|\asymp \frac{|\tau -w_\ve|}{t}\, . $$ Now
let $|\theta -\eta|< \ve /32$. Then, by the analogue of Lemma
\ref{lem1.3} (cf. \cite[Lemma 1]{ad}) for the conformal mapping
$\Phi$, we obtain that $$ |z-\xi |\le\frac 12\mb{ dist }(z,L), $$
and, consequently,
\begin{eqnarray*}
\left| \frac{z-\xi }{\z_1-\z_3}\right|&\asymp&
\left| \frac{z_L-z }{z_L-\z_1}\right| \,
\left| \frac{\xi -z }{z_L-z}\right| \\
&\p&
\frac{|\tau |-1}{t}\, \frac{|\tau -w_\ve|}{|\tau|-1}=
\frac{|\tau -w_\ve|}{t}\, .
\end{eqnarray*}
Hence, (\ref{n.7}) and  (\ref{4.12}) give
\beq \label{5.4}
|\tilde{r}(\tau )|\preceq \frac{|\tau -w_\ve |^2}{t^2},\quad |\tau
|=1+\ve , \, |\tau -w_\ve |\le c_2\, t. \eeq Relation (\ref{5.4})
remains true for $\tau$ such that $|\tau |>1+\ve , \, |\tau
-w_\ve|= c_2\, t$,  by the definition of the function
$\tilde{r}(\tau )$ and (\ref{4.10})--(\ref{4.12}).

Further, a direct computation shows that \beq \label{5.5}
|\tilde{r}(w)|\preceq t. \eeq Indeed, let us introduce the
auxiliary function $\tilde{R}(\tau)$, which we define to be the
harmonic extension of the function $$ \tilde{R}(\tau):=\left\{
\begin{array}{ll} |\tilde{r}(\tau)| , & \mb{if }\, |\tau |=1+\ve
,\, |\tau -w_\ve |\le c_2\, t,\\[2ex] c_3, & \mb{otherwise for
}|\tau |=1+\ve,
\end{array}\right.
$$ to $|\tau |\ge 1+\ve$. It is clear that we have, for
sufficiently large $c_3$, $$ |\tilde{r}(\tau )|\le \tilde{R}(\tau )
$$ on the boundary of the domain $$ \{ \tau :|\tau |>1+\ve ,\,
|\tau -w_\ve |<c_2\, t\}. $$ Therefore, by the maximum principle
for harmonic functions, the Poisson formula and (\ref{5.4}), we
obtain that
\begin{eqnarray*}
|\tilde{r}(w )|&\le& \tilde{R}(w ) =\tilde{R}(r e^{i\theta} )
\\
&=&\frac{1}{2\pi}\int_0^{2\pi}\tilde{R}((1+\ve )e^{i\eta} )\,
\frac{r^2-(1+\ve )^2}{r^2-2\, r\, (1+\ve )\, \mb{cos }(\theta -
\eta )+(1+\ve )^2}\, d\eta
\\
&\preceq & \ve \left( \frac{1}{t^2}\int_{\theta -c_2t}^{\theta +c_2t}
d\eta +\int_{\theta +c_2t}^{\theta +\pi}\frac{d\eta}{(\eta -\theta )^2}+
\int^{\theta -c_2t}_{\theta -\pi}\frac{d\eta}{(\eta -\theta )^2}
\right)
\preceq \frac{\ve}{t}=t.
\end{eqnarray*}

Comparing (\ref{4.10}), (\ref{4.11}),  (\ref{5.2}) and
(\ref{5.5}), we get the desired inequality (\ref{4.13}) by
(\ref{3.1}).

The same reasoning gives an analogue of (\ref{4.13}) for the case
$$ 1+\ve
<r<1+2\ve \, ,\quad
 \theta_2+\ve_1 t\le \theta \le \theta_2+2c_1t \, .
$$

Next, we assume that
\beq \label{n.6} 1<r=|w|<1+\ve ,\quad \z=\Psi
(w),\quad \z_\ve=\Psi (w_\ve). \eeq
Note that $L_\ve$ is convex
(cf. \cite[p. 47]{pom1}). Moreover, since $\Phi$ has a
quasiconformal extension to $\ov{\Bbb C}$, each $L_\ve$ is
$K$-quasiconformal with $K\ge 1$, independent of $\ve$. Therefore,
we have, by formula (\ref{n3.11}) for any $2|\z -\z_\ve |<s<\ve_2$
and any function $\kappa (z)$, continuous on $\ov{\mb{int }L_\ve}$
and harmonic in int $L_\ve$, that
\begin{eqnarray}
    |\kappa (\z ) - \kappa (\z_\ve )|
    &\preceq&
    \omega_{\z_\ve ,\kappa,L_\ve }(c_4\, |\z -\z_\ve |)
    \nonumber\\
\label{5.8}
&+&|\z -\z_\ve | \int^{s}_{|\z-\z_\ve |}
\frac{\omega_{\z_\ve ,\kappa,L_\ve }(c_4\, r)}
   {r^2}\, dr +
 \frac{|\z-\z_\ve |}{s}||\kappa ||_{L_\ve}\, ,
\end{eqnarray}
where $c_4>0$ is independent of $\z$ and $\ve$.

It is easy to prove (\ref{4.13}), if, in addition to (\ref{n.6}),
$\z_L\not\in J(2c_1t)\setminus J(\ve_1 t)$. Indeed, let now
$\kappa :=f,\, s:=\ve_3\, |\z_L-\z_L^*|$, where $\z_L^*:= \Psi
(e^{it}\Phi (\z_L))$ and the sufficiently small constant $\ve_3$
is chosen such that $\omega_{\z_\ve ,\kappa,L_\ve }(c_4\, s)=0$.
Therefore, we obtain (\ref{4.13}) by (\ref{5.8}), Lemma \ref{q}
and the obvious inequality $$ \left|
\frac{\z_L-\z_\ve}{\z_L-\z_L^*}\right| \preceq \frac{\ve }{t}=t\,
, $$ which follows from Lemma \ref{a}.

The situation is more complicated if, in addition to (\ref{n.6}),
$\z_L\in J(2c_1t)\setminus J(\ve_1 t)$. For definiteness, let
$\z_L\in L(\z_3,\z_1)$. In this case, we represent the function
$g$ in the form of (\ref{4.10}) with $\xi :=\z_\ve$, and set
$\kappa (z):=r (z,\xi )$ (i.e. $\kappa (z)$ is the harmonic
extension of $r (z,\xi )$ from $L_\ve$ to int$L_\ve$),
$\,s:=\ve_4\, |\z_1-\z_3|$, where $\ve_4$ is chosen to be so small
that the function $\kappa (z)$ satisfies (\ref{4.12}) for $z\in
l(s)$. Since  $$ |\kappa (z)|\preceq 1,\quad z\in \gamma (s), $$
by (\ref{4.10}) and (\ref{4.11}), we have on setting
$\delta:=|\z-\z_\ve|$ that
\begin{eqnarray}
|r(\z ,\xi )|&\preceq& \frac{\delta^2}{s^2}+\frac{\delta}{s^2}\int_{\delta}^sdr
+\frac{\delta}{s}
\nonumber\\
\label{5.9}
 & \asymp& \frac{\delta}{s}\p \left|
 \frac{\z_L-\z_\ve}{\z_1-\z_3}\right|,
\end{eqnarray}
by (\ref{4.12}) and (\ref{5.8}). Comparing (\ref{5.9}),
(\ref{4.10}), (\ref{4.11}) and applying Lemma \ref{a}, we get $$
|f(\z )-f(\z_\ve)|\preceq \left|
\frac{\z_L-\z_\ve}{\z_1-\z_3}\right| \asymp \left|
\frac{\z_L-\z_\ve}{\z_L-\z_1}\right| \preceq \frac{\ve}{t}=t. $$

\hfill\rbox

\bigskip

V. V. Andrievskii, Mathematisch-Geographische Fakultaet, Katholische Universitaet Eichstaett,
D-85071 Eichstaett, Germany.

I. E. Pritsker, Department of Mathematics, Case Western Reserve University,
10900 Euclid Avenue, Cleveland, Ohio 44106-7058, U.S.A.

R. S. Varga, Institute for Computational Mathematics, Department of Mathematics
and Computer Science, Kent State University, Kent, OH 44242, U.S.A.


\begin{thebibliography}{99}

\bibitem{ahl} L.V. Ahlfors (1966), Lectures on Quasiconformal Mappings.
               Princeton, N.J.: Van Nostrand.

\bibitem{andbeldzj} V.V. Andrievskii, V.I. Belyi, V.K.
Dzyadyk (1995),
Conformal Invariants in Constructive Theory of
Functions of Complex Variable. Atlanta, Georgia: World Federation Publisher.

\bibitem{andbla} V.V. Andrievskii, H.-P. Blatt (1997), {\it
Erd\H{o}s-Tur\'an type
theorems on piecewise smooth curves and arcs.} J. Approx. Theory, {\bf 88},
109-134.

\bibitem{ad} V.V. Andrievskii, H.-P. Blatt (1997), {\it A
discrepancy theorem on quasiconformal curves.} Constr. Approx.,
{\bf 13}, 363-379.

\bibitem{andbla1} V.V. Andrievskii, H.-P. Blatt (1999),
{\em  Erd\H{o}s-Tur\'an type theorems on quasiconformal curves and
arcs.} J. Approx. Theory, to appear.

\bibitem{belmik} V.I. Belyi, V.M. Mikljukov (1974), {\it Some properties
of conformal and quasiconformal mappings and direct theorems of the constructive
theory of functions.} Izv. Akad. Nauk SSSR Ser. Mat., {\bf 38},
1343--1361. (Russian)

\bibitem{blagro} H.-P. Blatt,  R. Grothmann (1991), {\it Erd\H{o}s-Tur\'an
                theorems on a system of Jordan curves and arcs}. Constr.
                Approx., {\bf 7}, 19-47.

\bibitem{blamha} H.-P. Blatt,  H. Mhaskar (1993), {\it A general discrepancy
                 theorem.} Ark. Mat., {\bf 31}, 219-246.

\bibitem{eiesta}  M. Eiermann, H. Stahl (1994),
{\em Zeros of orthogonal polynomials on regular N-gons}, in
``Linear and Complex Analysis Problem Book 3" (V.P. Havin and N.K.
Nikolski, eds.), vol. 2, Lecture Notes in Math. {\bf 1574},
187-189.

\bibitem{erdtur} P. Erd\H{o}s, P. Tur\'an (1950), {\it On the distribution
                of roots of polynomials.} Annals of Math., {\bf 51}, 105-119.

\bibitem{gan} T. Ganelius (1957), {\it Some applications of a lemma
on Fourier series.} Academie Serbe des Sciences Publications de
l'Institut Matematique (Belgrade), {\bf 11}, 9-18.

\bibitem{hal}
K. Haliste (1967), {\em Estimates of harmonic measures.} Ark.
Math., {\bf 6}, 1--31.

\bibitem{her}
J. Hersch (1955), {\em Longueurs extr\'emales et theory des
founctions.} Comm. Math. Helv., {\bf 29}, 301--337.

\bibitem{koepom} T. K\H{o}vari, Ch Pommerenke (1967),
{\em On Faber polynomials and Faber expansions.} Mathematische
Zeitschrift, {\bf 99}, 193--206.

\bibitem{leh} O. Lehto (1987), Univalent Functions and Teichmller
              Spaces. New York: Springer-Verlag.

\bibitem{lehvir}
O. Lehto, K. I. Virtanen (1973),  Quasiconformal Mappings in the
Plane, 2nd ed. New York: Springer-Verlag.

\bibitem{pom1}
Ch. Pommerenke (1975), Univalent Functions.
                G\"{o}ttingen: Vandenhoeck and Ruprecht.

\bibitem{pom2}
Ch. Pommerenke (1992),  Boundary Behaviour of Conformal Maps.
                New York: Springer-Verlag.

\bibitem{saftot} E. Saff, V. Totik (1997), Logarithmic Potentials
with External Fields. New York/Berlin: Springer-Verlag.

\bibitem{smileb} V.I. Smirnov, N.A. Lebedev (1968), Functions of a Complex
Variable. Constructive Theory. Cambridge: Massachusetts Institute of
Technology.

\bibitem{statot}  H. Stahl, V. Totik (1992), General Orthogonal
Polynomials. Vol. 43 of Encyclopedia of Mathematics. Cambridge
University Press, New York.

\bibitem{ste} E.M. Stein (1970), Singular Integrals and Differentiability
Properties of Functions. Princeton, NJ: Princeton University Press.

\bibitem{tot}  V. Totik (1993), {\it Distribution of simple zeros of
                 polynomials}.
                  Acta Math., {\bf 170}, 1--28.

\bibitem{ull} J. L. Ullman (1972),
{\em The location of the zeros of the derivatives of Faber
polynomials.} Proc. Amer. Math. Soc., {\bf 34}, 422--424.

\end{thebibliography}
\end{document}